\documentclass{article}
\usepackage{amssymb}

\newcommand{\R}{{\mathbb R}}
\newcommand{\N}{{\mathbb N}}

\newcommand{\C}{{\mathbb C}}

\newcommand{\A}{{\cal A}}

\newcommand{\cL}{{\cal L}}

\newcommand{\M}{{\cal M}}

\newcommand{\cH}{{\cal H}}

\newcommand{\gG}{\Gamma}
\newcommand{\gve}{\varepsilon}
\newcommand{\gl}{\lambda}

\newcommand{\gt}{\theta}
\newcommand{\gT}{\Theta}

\begin{document}
\noindent{\bf \large Construction de Triplets Spectraux \`a Partir de Modules 
de Fredholm}

\vspace{1ex}
\noindent{\bf Elmar {SCHROHE}, Markus {WALZE}%
\footnote{support\'e par la bourse `TMR Marie Curie Research Training Grants' 
 contract No. ERBFMBICT$961518$} et Jan-Martin {WARZECHA}}

\vspace{1ex}
\noindent{\small \bf Institut f\"ur Mathematik, Universit\"at Potsdam, 
14415 Potsdam, Allemagne}\\
{\small \bf IHES, 35, route de Chartres, 91440 Bures-sur-Yvette, France}\\
{\small \bf Institut f\"ur Physik, Johannes Gutenberg-Universit\"at Mainz, 
55099 Mainz, Allemagne}

\vspace{2ex}
{\small
{\bf \it R\'esum\'e} -- Soit $(\A,\cH, F)$ un module de Fredholm $p$-sommable,
o\`u l'alg\`ebre $\A= \C \Gamma$ est engendr\'ee par un groupe discret $\Gamma$
d'\'el\'ements unitaires de ${\cal L}(\cH)$ qui est de croissance 
polynomiale $r$.  On construit alors un triplet spectral $(\A,\cH,D)$ de 
sommabilit\'e $q$ pour tout $q>p+r+1$ avec $F={\rm sign} \,D$. 
Dans le cas o\`u $(\A,\cH,F)$ est 
$(p,\infty)$-sommable on obtient la $(q,\infty)$-sommabilit\'e de $(\A,\cH,D)$ pour tout
$q>p+r+1$.}

\vspace{2ex}

\vspace{1ex}
{\small
{\bf \it Abstract} -- Let $(\A,\cH,F)$ be a $p$-summable Fredholm module where
the algebra $\A=\C \Gamma$ is generated by a discrete group of unitaries in 
${\cal L}(H)$ which is of polynomial growth $r$. Then we construct a 
spectral triple  $(\A,\cH,D)$ with $F= {\rm sign}\,D$ which is 
$q$-summable for each $q>p+r+1$. In 
case $(\A,\cH,F)$ is $(p,\infty)$-summable we obtain 
$(q,\infty)$-summability of $(\A,\cH,D)$ for each $q>p+r+1$.}

\vspace{2ex}
\noindent{\bf Abridged English Version}

In [\ref{Connes91},~Th\'eor\`eme 3] A.~Connes showed the following theorem
which we quote as stated in [\ref{Connes94},~IV. Theorem 4]:\medskip

{\bf Theorem}. {\sl Let $A$ be a $C^*$-algebra, $(\cH,F)$ a Fredholm module
over $A$ and $\A \subseteq A$ a countably generated subalgebra such that,
 for each $ a \in \A$,
\begin{equation}
[F,a] \in {\rm Li}^{1/2} (\cH).
\end{equation}
Then there exists a self-adjoint unbounded operator $D$ in $\cH$ such that
\hspace*{2em}\begin{itemize}
\item[] {\rm(D1)}\quad${\rm Sign}\, D = F$,\\[-3ex]
\item[] {\rm(D2)}\quad$[D,a]$ is bounded for any $a \in \A$,\\[-3ex]
\item[] {\rm(D3)}\quad${\rm Trace}\, (e^{-D^2}) < \infty$.
\end{itemize}}\bigskip

Here, ${\rm Li}^{1/2} (\cH)$ is the ideal of all compact operators on $\cH$
whose singular values $\{ \mu_n \}^\infty_{n=0}$ satisfy $\mu_n = O((\log
n)^{-1/2})$ as $n \to \infty$. Similarly ${\rm Li} (\cH)$ is characterized by 
the property $\mu_n = O((\log n)^{-1})$. The algebra $B = \{ T \in A: [F,T] 
\in {\rm Li}^{1/2} (\cH) \}$ is symmetric and stable under holomorphic 
functional calculus so that one can enlarge $\A$ and assume that it is 
generated by a countable discrete group of unitaries, $\Gamma$. 

For simplicity we shall start from this point of view. We suppose that
$\A=\C \Gamma$, where $\Gamma$ is a discrete group of unitaries in $\cL(\cH)$ 
which is of polynomial growth $r$. Moreover, we assume that $(\A,\cH,F)$ is a 
$p$-summable (respectively $(p,\infty)$-summable) Fredholm module over $\A$, 
i.~e.,
\begin{equation}
[F,a] \in \cL^p (\cH) \qquad {\rm (respectively}~ \cL^{(p,\infty)} (\cH)) \quad
{\rm for~all~} a \in \A .
\end{equation}
We show that we can find an unbounded self-adjoint operator $D$ satisfying 
(D1), (D2), and, for each $q>p+r+1$, 
\begin{itemize}
\item[](D3')\ \ \ \ \ $(1+D^2)^{-q/2} \in \cL^1 (\cH)$ or\\[-3ex]
\item[](D3'')\ \ \ \ \ $(1+D^2)^{-q/2} \in \cL^{(1,\infty)} (\cH)$, 
respectively.
\end{itemize}
The triple $(\A,\cH,D)$ is called a $q$-summable spectral triple (and 
$(q,\infty)$-summable in case (D3'')). Spectral triples are sometimes  
called {\em unbounded Fredholm modules}.

Our proof essentially follows the original idea of Connes; additional 
ingredients are 
\begin{itemize}
\item a characterization of those selfadjoint operators $|D|$ for which
 which both $[|D|,a]$ and $[F,a]|D|$ are bounded (Proposition 3),  
 \\[-3ex]
\item a different function for the nonlinear transformation (see below), and
\item
 a theorem of Rotfel'd to estimate the singular values of $D$. 
\end{itemize}

If $\gG$ is  finitely generated then boundedness of $[D,a]$ will even hold for 
all $a\in C_1(\gG)$, cf.~[\ref{Connes91}, D\'efinition 5].
\bigskip

Nous rappelons qu'une  fonction longueur $L$ sur un groupe discret $\Gamma$ est
une application $L: \Gamma \to [0,\infty)$ qui v\'erifie, pour tous 
$g,h \in \Gamma$,
$$
L(gh) \leq L(g) + L(h), \quad L(g^{-1}) = L(g), \quad L(1) = 0 .
$$
Si, de plus, le cardinal de l'ensemble $B_k = \{ g \in \Gamma: L(g) \leq
k \}, k \in \N_0$, est $O((1+k)^r)$, alors $\Gamma$ est dit de croissance
polynomiale d'ordre $r$.

Soit $p>1$ un nombre r\'eel. Un module de Fredholm $p$-sommable (respectivement
$(p,\infty)$-sommable), est un triplet $(\A,\cH,F)$ o\`u $\A$ est une alg\`ebre 
unif\`ere, $\cH$ est un espace de Hilbert avec une repr\'esentation 
$\pi: \A \to \cL(\cH)$, et $F=F^* \in \cL(\cH)$ est un op\'erateur 
satisfaisant $F^2=I$ et, pour tout $a \in \A$,
\begin{equation}
[F,a] \in \cL^p (\cH) \quad ({\rm resp.} \quad [F,a] \in \cL^{(p,\infty)} (\cH)).
\end{equation}
On parle d'un module de Fredholm $\theta$-sommable, si (3) est remplac\'e par 
(1).

Un triplet spectral $p$-sommable (respectivement $(p,\infty)$-sommable)
$(\A,\cH,D)$ est constitu\'e d'une alg\`ebre unif\`ere $\A$, repr\'esent\'e 
sur un espace de Hilbert $\cH$, et d'un op\'erateur autoadjoint $D$ \`a 
r\'esolvante compacte tel que, pour tout $a \in \A$,
$$
[D,a] \in \cL(\cH)
$$
et
\begin{equation}
(1+D^2)^{-p/2} \in \cL^1 (\cH) \quad {\rm (resp.~}
(1+D^2)^{-p/2} \in \cL^{(1,\infty)} (\cH)) .
\end{equation}
Si la condition (4) est remplac\'ee par
$$
e^{-D^2} \in \cL^1 (\cH)
$$
on parle de $\theta$-sommabilit\'e.
\bigskip

{\bf 1.~Proposition.} {\sl Soit $(\A,\cH,D)$ un triplet spectral et $F = {\rm
sign}\, D$. Alors $(\A,\cH,F)$ est un module de Fredholm. Si $(\A,\cH,D)$
est $p$-sommable {\rm(}resp.\ $(p,\infty)$-sommable, resp.\
$\theta$-sommable{\rm)} alors  $(\A,\cH,F)$ est $p$-sommable
{\rm(}resp.\ $(p,\infty)$-sommable, resp.\ $\theta$-sommable{\rm)}.}
\bigskip

{\em D\'emonstration.} Sans restreindre le cas g\'en\'eral on peut supposer 
que $D$ est inversible. 
En utilisant la formule $T^{-1/2} = \frac1\pi \int_0^\infty \gl^{-1/2} (\gl+T
)^{-1}\, d\gl$ on montre que
\begin{eqnarray*} 
\lefteqn{[F,a] =  \frac1\pi \int_0^\infty \gl^{-1/2} D[(\gl+D^2)^{-1},a]\, d\gl
+[D,a]|D|^{-1}}\\
&=& \frac1\pi \int_0^\infty\left( \gl^{1/2} (\gl+D^2)^{-1}[D,a] (\gl+D^2)^{-1}
-\gl^{-1/2} D (\gl+D^2)^{-1}[D,a]  (\gl+D^2)^{-1}D\right)\, d\gl.
\end{eqnarray*}
Puisque tout \'el\'ement de l'alg\`ebre sym\'etrique $\A$ est la somme d'un 
\'el\'ement autoadjoint et d'un \'el\'ement anti-autoadjoint on peut supposer 
que $[D,a]$ est autoadjoint.

En utilisant les in\'egalit\'es $\, -\|[D,a]\| \le [D,a]\le \|[D,a]\|$ et la formule pour 
l'inverse de la racine, on obtient 
$$-\|[D,a]\| |D|^{-1} \le [F,a] \le \|[D,a]\||D|^{-1}.$$ 
De l\`a d\'ecoule la proposition. $\square$
\bigskip

R\'eciproquement, Connes  [\ref{Connes89}] et Voiculescu [\ref{Voiculescu90}]
ont montr\'e qu'il existe des obstructions \`a l'existence de triplets 
spectraux $(\A,\cH,D)$ de sommabilit\'e finie associ\'es \`a une alg\`ebre 
donn\'ee, $\A$. En particulier, si $\Gamma$ est un groupe discret non-moyennable,
alors il n'existe pas de triplet spectral de sommabilit\'e finie associ\'e \`a 
une sous-alg\`ebre dense de $C^*_{\rm red}(\Gamma)$. Il y a m\^eme des 
obstructions pour un groupe resoluble de croissance exponentielle (donc 
moyennable). D'autre part, si $\Gamma$ est de type fini et de croissance 
polynomiale $r$ alors on peut trouver un triplet $(\A,\cH,D)$ pour
$\A=\C\Gamma, \cH = l^2 (\Gamma)$. En effet $D$ est l'op\'erateur positif 
agissant par multiplication par la fonction longueur $L$, et l'on a 
$(1+D^2)^{-p/2} \in \cL^1 (\cH)$ pour tout $p>r+1$.\medskip

Le r\'esultat principal de cette Note est le th\'eor\`eme suivant.
\bigskip

{\bf 2.~Th\'eor\`eme.} {\sl Soit $\Gamma$ un groupe discret de croissance
polynomiale $r$, repr\'esent\'e par des \'el\'ements unitaires sur un
espace de Hilbert $\cH$. Si $(\A,\cH,F)$ est un module de Fredholm
$p$-sommable {\rm(}resp.\ $(p,\infty)$-sommable{\rm)} avec $\A=\C\Gamma$, 
alors il existe un triplet spectral $(\A,\cH,D)$ de sommabilit\'e $q$ 
{\rm(}resp.\ $(q,\infty)${\rm)}  pour tout $q>p+r+1$ tel que $D=F|D|$.

Pour le cas o\`u $\gG$ est le groupe abelien libre engendr\'e par $r$ 
\'el\'ements il est possible d'obtenir la sommabilit\'e $q$ 
resp.~$(q,\infty )$ pour tout $q>p+r$}.      
\bigskip

Il est essentiel, ici, que $D$ soit li\'e \`a $F$ par $D=F|D|$, sinon on
pourrait obtenir une meilleure sommabilit\'e selon les r\'esultats mentionn\'es
ci-dessus.

La d\'emonstration suit l'id\'ee de Connes. On choisit des g\'en\'erateurs 
$u^1,u^2,\ldots$ de $\Gamma$ et l'on introduit la ``m\'etrique quantique''
\begin{equation}
G= \sum_k c_k [F,u^k]^* [F,u^k]
\end{equation}
avec $c_k \in \R_+$ et  $c_k\le 2^{-k}\|[F,u^k]\|^{-2}$ avec la norme dans 
$\cL^p(\cH)$ (resp.~ dans $\cL^{(p,\infty)}(\cH)$). 

On construira $|D|$ a partir de $G$. En effet, on utilisera l'op\'erateur 
$$
\Theta (G)=f^{-1}{\cal M}f(G)    
$$
avec une moyennisation ${\cal M}$ et une transformation non lin\'eaire $f$;
ensuite on posera $D=F|D|$. Pour v\'erifier que le commutateur $[D,a]$ est 
born\'e pour tout $a\in\C\gG$, il est donc suffisant que $[|D|,a]$ et
$[F,a]|D|$ soient born\'es. Nous effectuons une observation importante :
\bigskip

{\bf 3.~Proposition.} {\sl Soit $T$ un op\'erateur  autoadjoint inversible 
{\rm(}\'eventuellement non-born\'e{\rm)}. Alors les conditions suivantes sont 
\'equivalentes :
\renewcommand{\labelenumi}{{\rm(\roman{enumi})}}
\begin{enumerate}
\item $[T,a]$ et $[F,a] T$ sont born\'es pour tout $a \in \C\Gamma$.
\item Pour chaque op\'erateur unitaire $u \in \Gamma$ il existe une
      constante $C_u>0$ telle que
      $$
      T^{-1} (I-C_u T^{-1}) \leq u T^{-1} u^* \leq T^{-1} (I+C_u
      T^{-1}) .
      $$
      De plus il existe $\gl>0$ tel que $T^{-2} \geq \gl G$.
\end{enumerate}}
Dans un premier temps, on choisit une fonction $f \in C^\infty [0,\gve), \gve
>0$, avec $f(0)=0$, $f$ strictement croissante et telle que l'inverse
$f^{-1}$ soit une fonction concave et croissante d'op\'erateurs.

Fixons aussi la fonction de poids $\rho: \Gamma \to (0,1)$ donn\'ee par
$\rho(u) = \exp (-(1+L(u)))$ et l'op\'erateur de moyennisation, $\M$, defini
par $\M(T) = \sum_{u \in \Gamma} \rho(u) u T u^{*}$ pour $T\in\cL(\cH)$. A 
l'aide de $f$ et $\M$ nous introduisons l'op\'erateur $\Theta(G) = 
f^{-1}{\cal M}f(G).$

Il est clair que $\Theta(G)$ est  compact et positif. Supposons de plus que 
$\Theta(G)$ est injectif. Cela nous permet de poser 
$$|D| = \Theta(G)^{-1/2}+F\Theta(G)^{-1/2}F$$ 
et $D = F|D|$ (notons que $[F,|D|] = 0$).
 
En utilisant la proposition 3 on peut montrer que $D$ a les propri\'et\'es 
d\'esir\'ees pourvu que les trois conditions suivantes soient v\'erifi\'ees:
\begin{itemize}
\item[] (T1)\quad $\Theta(G) \geq \gl G$ pour un $\gl>0$;
\item[] (T2)\quad $\Theta(G)^{1/2} (I- C_u \Theta(G)^{1/2}) \leq u
            \Theta(G)^{1/2} u^* \leq \Theta(G)^{1/2} (I+C_u
            \Theta(G)^{1/2})$;
\item[] (T3)\quad $\Theta(G) \in \cL^{q/2} (\cH)$ (resp.\ $\Theta(G) \in
            \cL^{(q/2,\infty)} (\cH)$), $q>p+r+1$.
\end{itemize}
La propri\'et\'e (T1) est ais\'ement d\'emontr\'ee. En revanche, le choix de 
$f$ est essentiel pour montrer que (T2) et (T3) sont verifi\'ees.
Naturellement on essaiera d'atteindre une valeur de $q$ proche de $p$ dans la
relation (T3). D'apr\`es la proposition 1 on aura $q\ge p$. Cependant, il n'est 
pas \'evident de d\'eterminer la valeur minimale de $q$, et il faut
bien choisir $f$ pour obtenir $q>p+r+1$. Nous utilisons la fonction dont 
l'inverse est donn\'e par
$$
f^{-1} (t) = \left( {\rm arcosh}\, \frac{1}{t} \right)^{-2} = \left( \log \left(
\frac{1}{t} \left( 1 - \sqrt{1-t^2} \right) \right) \right)^{-2} .
$$
On remarque que
\begin{equation}
f^{-1} \sim (\log t)^{-2}, \qquad t>0 \quad {\rm petit},
\end{equation}
par cons\'equent, les fonctions $f^{-1}$ et $f$ sont croissantes pr\`es de
zero. Ensuite on utilise la caract\'erisation de L\"owner [\ref{Loewner}] pour
verifier que $f^{-1}$ est une fonction croissante d'op\'erateurs, 
c'est-\`a-dire qu'on d\'emontre que $f^{-1} : [0,\gve) \to \R$ s'\'etend \`a une 
fonction analytique dans le demi-plan $\{ z \in \C: {\rm Im}\, z>0\}$. La 
relation (6) et le fait que $f^{-1}$ est une fonction croissante 
d'op\'erateurs entra\^{\i}nent (T2).

Finalement, on designe par $\{ \mu_0 (T), \mu_1 (T), \ldots \}$ les
valeurs singuli\`eres d'un op\'erateur compact, $T$.

Selon un th\'eor\`eme de Rotfel'd [\ref{Rotfeld}], la concavit\'e de $f^{-1}$ 
implique l'in\'egalit\'e
$${\sum^N_{m=0} \mu_m (\Theta (G))^q \leq \sum^N_{m=0} \sum_{u \in
\Gamma} f^{-1} (\mu_m (\rho(u) u f(G) u^{-1}))^q}  
=\sum^N_{m=0} \sum_{u \in \Gamma} f^{-1} (\rho(u) f(\mu_m(G)))^q 
$$
pour tout $N\in\N$. En utilisant la relation $\rho(u) = \exp(- (1+L(u)))$, la 
croissance polynomiale de $\Gamma$ et le fait que $f^{-1} (t) \sim 
(\log t)^{-2}$, on en d\'eduit que
$$
\sum_m\mu_m (\Theta(G))^q \leq C \sum_m\mu_m (G)^{q-\frac{1}{2}(r+1)} .
$$
De cette estimation d\'ecoule (T3) et donc l'\'enonc\'e du th\'eor\`eme.

Il reste \`a consid\'erer le cas o\`u $\gT(G)$ a un noyau non trivial
d\^u \`a ${\rm ker}\, G$. On pose 
$$\cH_0 = {\rm ker}\,\gT(G) = \bigcap_{u,v\in\gG}{\rm ker}([F,u]v)=
\bigcap_{u\in\gG}{\rm ker}([F,u]).
$$ 
Alors on peut \'ecrire $\cH = \cH_0 \oplus \cH_1$ o\`u $\cH_0$ est invariant 
par $\A$ et $F$, o\`u $F|_{\cH_0}$ commute avec tout $a \in \A$ et o\`u 
$\Theta(G)$ est injectif sur $\cH_1$. La construction de $D$ sur $\cH_0$ sera 
facile: au lieu de $G$ on choisira un \'el\'ement arbitraire $G_0\in 
\cL^{p/2}(\cH_0)$ (resp.~${\cL}^{(p/2,\infty)}$) qui est strictement positif, 
ensuite on appliquera la m\^eme consid\'eration que sur $\cH_1$.
\medskip

{\bf 4.~Remarque.} Une construction similaire peut \^etre utilis\'ee dans le
cas $\gt$-sommable.

{\bf 5.~Remarque.} Il est \'evident que $[D,a]$ sera born\'e pour tout $a$ 
dans la compl\'etion de $\C\gG$ pour la norme 
$$\|a\|^{\sim} = \|a\|_{{\cal L}(H)}+\|[D,a]\|_{{\cal L}(H)}.$$

Pour $\gG$ de type fini, $[D,b]$ 
sera born\'e pour tout $b\in C_1(\gG)$, cf.~[\ref{Connes91}, D\'efinition 5],
puisque, pour tout $a\in\C\gG$, la norme $\|[D,a]\|$ 
peut \^etre estim\'ee par $c\|a\|'$, o\`u $\|\cdot\|'$ est la norme 
utilis\'ee dans la construction de $C_1(\gG)$ et $c$ est une constante 
universelle.
\bigskip

\small
{\sc R\'ef\'erences bibliographiques} 
\renewcommand{\labelenumi}{[\arabic{enumi}]}
\begin{enumerate}

\item\label{Connes89}
A.~Connes. 
\newblock Compact metric spaces, {F}redholm modules, and hyperfiniteness.
\newblock {\em Ergod. Th. \& Dynam. Sys.} {\bf 9} (1989), 207--220.

\item \label{Connes91}
A.~Connes.
\newblock Caract\`eres de repr\'esentations $\theta$-sommables des groupes
  discrets.
\newblock {\em C. R. Acad. Sci. Paris} {\bf 312}, s\'erie I  (1991), 661--666.

\item \label{Connes94}
A.~Connes.
\newblock {\em Noncommutative Geometry}.
\newblock Academic Press, New York, London, Tokyo 1994.

\item \label{Loewner}
K.~L\"owner.
\newblock {\"U}ber monotone {M}atrixfunktionen.
\newblock {\em Math. Z.} {\bf 38} (1934), 177--216.

\item \label{Rotfeld}
S.~Yu. Rotfel'd.
\newblock The singular values of the sum of completely continuous operators.
\newblock In: Sh. Birman (Ed.), {\em Spectral Theory}. Topics
  in Mathematical Physics, vol. 3,  
  Consultants Bureau, New York, London 1963, pp.~73 -- 78.

\item \label{SchroheWalzeWarzecha}
E.~Schrohe, M.~Walze, and J.-M. Warzecha.
\newblock Lifting bounded {F}redholm modules to spectral triples.
\newblock En pr\'eparation.

\item \label{Voiculescu90}
D.~Voiculescu. 
\newblock On the existence of quasicentral approximate units relative to normed
ideals. Part I.
\newblock {\em J. Funct. Anal.} {\bf 91} (1990), 1--36.
\end{enumerate}
\bigskip

\noindent{\em schrohe@mpg-ana.uni-potsdam.de\\
walze@ihes.fr\\
warzecha@thep.physik.uni-mainz.de}
\end{document}